\documentclass[a4paper, reqno]{amsart}

\usepackage[utf8]{inputenc}
\usepackage{mathtools}
\usepackage{graphicx}
\usepackage{array}
\usepackage{url}
\usepackage{tikz}
\newcommand{\R}{\mathbb{R}}

\DeclareMathOperator{\vol}{vol}

\date{11 March 2019}

\title{A counterexample to a conjecture of Larman and Rogers on sets
  avoiding distance~1}

\author{Fernando Mário de Oliveira Filho}
\address{F.M. de Oliveira Filho, Delft Institute of Applied
  Mathematics, Delft University of Technology, Van Mourik Broekmanweg
  6, 2628 XE Delft, The Netherlands.}
\email{fmario@gmail.com}

\author{Frank Vallentin}
\address{F.~Vallentin, Mathematisches Institut, Universit\"at zu
  K\"oln, Weyertal~86--90, 50931 K\"oln, Germany.}
\email{frank.vallentin@uni-koeln.de}

\thanks{The second author was partially supported by the SFB/TRR 191
  ``Symplectic Structures in Geometry, Algebra and Dynamics'', funded
  by the DFG, and the European Union's Horizon 2020 research and
  innovation programme under the Marie Sk\l{}odowska-Curie agreement
  number~764759}

\subjclass[2010]{52C10, 51K99}

\begin{document}

\begin{abstract}
  For each~$n \geq 2$ we construct a measurable subset of the unit
  ball in~$\R^n$ that does not contain pairs of points at distance~1
  and whose volume is greater than~$(1/2)^n$ times the volume of the
  unit ball. This disproves a conjecture of Larman and Rogers
  from~1972.
\end{abstract}

\maketitle
\markboth{F.M. de Oliveira Filho and F. Vallentin}{A counterexample to
  a conjecture of Larman and Rogers}

%%%%%%%%%%%%%%%%%%%%%%%%%%%%%%%%%%%%%%%%%%%%%%%%%%%%%%%%%%%%%%%%%%%%%%

Larman and Rogers~\cite[Conjecture 1]{LarmanR1972} conjectured:
``Suppose that the distance~$1$ is not realized in a closed subset $S$
of a spherical ball $B$ of radius $1$. Then the Lebesgue measure of
$S$ is less than $(1/2)^n$ times the Lebesgue measure of $B$.'' Since
an open ball of radius~$1/2$ does not have pairs of points at
distance~1, this bound would be tight. Croft, Falconer, and
Guy~\cite[page 178]{CroftFG1991} comment on the planar case of this
conjecture, saying about the optimal set~$S$ (they disregard the
requirement that it has to be closed): ``Surely it must be a disk of
radius~$1/2$, but this seems hard to prove.''

The following simple construction shows that the conjecture is wrong
in all dimensions~$n \geq 2$. Let~$e_1 = (1, 0, \ldots, 0) \in \R^n$
and write~$a = (1/6) (1 + \sqrt{10})$. Let
\[
  T_n = \{\, x \in \R^n : x_1 > 1/2,\ \|x-a e_1\| < 1/2,\ \|x\| < 1\,\}
\]
and set~$S_n = T_n \cup -T_n$. The set~$T_n$ is the intersection of
the open ball of radius~$1/2$ centered at~$a e_1$ with the open ball
of radius~1 centered at the origin and the open halfspace~$x_1 >
1/2$. It is easy to see that~$T_n$ does not contain pairs of points at
distance~1, and hence neither does~$S_n$. This is the counterexample
to the conjecture\footnote{Since Larman and Rogers ask for a closed
  set, take closed inner approximations of~$S_n$.}, as is shown below;
see also Figure~\ref{fig:mosquito}.

\begin{figure}[t]
  \begin{center}
    \includegraphics{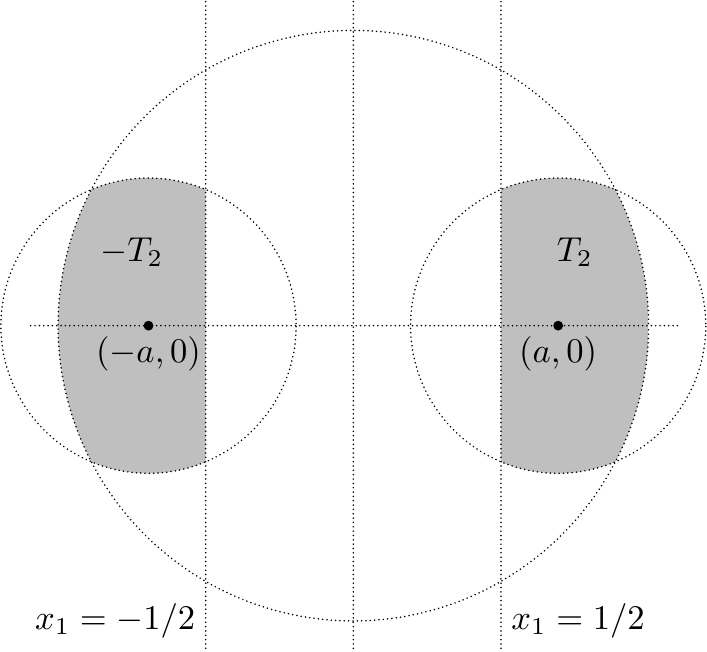}
  \end{center}

  \caption{The set~$S_2$. The parameter~$a$ is chosen so that~$a e_1$
    is equidistant to the hyperplane~$x_1 = 1/2$ and the hyperplane
    that contains the intersection between the spheres~$\|x\| = 1$ and
    $\|x - a e_1\| = 1/2$; this choice for~$a$ maximizes the volume
    of~$T_n$.}\label{fig:mosquito}
\end{figure}

The volume of the unit ball~$B_n$ is~$v_n = \pi^{n/2} / \Gamma(1 +
n/2)$. So
\[
  \vol T_n = \int_{1/2 - a}^{a - 1/2} v_{n-1} (1/4 - x^2)^{(n-1)/2}\, dx
  + \int_{2a - 1/2}^1 v_{n-1}(1 - x^2)^{(n-1)/2}\, dx.
\]
This gives~$\vol S_2 / \vol B_2 = 0.2848\ldots$
and~$\vol S_3 / \vol B_3 = 0.1563\ldots$.

Actually, for~$n \geq 3$ it suffices to use the lower bound
\[
  \vol T_n \geq \int_{1/2 - a}^{a - 1/2} v_{n-1} (1/4 -
  x^2)^{(n-1)/2}\, dx
\]
to disprove the conjecture. It is known that the volume of the unit
ball is concentrated around the equator\footnote{See e.g.\ Theorem~2.7
  in Blum, Hopcroft, and Kannan~\cite{BlumHK2018};
  Ball~\cite{Ball1997} and Matoušek~\cite{Matousek2002} present
  analogous results for the sphere instead of the ball.}:
if~$n \geq 3$ and~$c \geq 1$, then
\[
  \frac{\vol\{\, x \in B_n : |x_1| \leq c / \sqrt{n-1}\,\}}{\vol B_n}
  \geq 1 - (2 / c) e^{-c^2 / 2}.
\]
From this inequality one gets immediately the asymptotic relation 
\begin{equation}
  \label{eq:asymp}
  \frac{\vol S_n}{\vol B_n} = (2 - o(1)) (1/2)^n > (1/2)^n.
\end{equation}
In other words, in high dimension it is almost possible to fit two
balls of radius~$1/2$ inside the unit ball instead of only one like
conjectured by Larman and Rogers.  Choosing an appropriate
constant~$c$, one shows that~$\vol S_n / \vol B_n > (1/2)^n$ for
all~$n \geq 15$; the remaining cases can be checked directly.

One of the original motivations of Larman and Rogers in proposing
their conjecture is that it is related to another, still open
conjecture of L.~Moser, also popularized by Erd\H{o}s, on the global
behavior of sets avoiding distance~1. This conjecture states that the
upper density\footnote{The \textit{upper density} of a
  Lebesgue-measurable set~$X \subseteq \R^n$ is
  \[
    \sup_{p\in\R^n}\limsup_{T\to\infty} \frac{\vol(X \cap (p +
      [-T, T]^n))}{\vol [-T,T]^n}.
  \]
  Intuitively, it is the fraction of space covered by~$X$.}  of any
measurable subset of~$\R^2$ containing no pair of points at distance~1
is less than~$1/4$; Larman and Rogers' conjecture would imply that any
such subset of~$\R^2$ has upper density at most~$1/4$. Note that
Larman and Rogers' conjecture, if it were true, still would not imply
Moser's conjecture. Indeed, Larman and Rogers' conjecture says that a
\textit{closed} subset of the unit disk that avoids distance~1 has
area \textit{less than}~$1/4$ times the area of the unit disk; this in
turn implies that the area of a measurable subset of the unit disk
that avoids distance~1 is \textit{at most}~$1/4$ times the area of the
unit disk.

The construction of~$S_n$ shows that the local behavior resembles the
double cap conjecture~\cite[Conjecture~2.8]{Kalai2015}, which states
that the union of two antipodal spherical caps of radius~$\pi/4$ each
is a maximum-area subset of the unit sphere having no pairs of
orthogonal vectors; see DeCorte, Oliveira, and
Vallentin~\cite{DeCorteOV2018} for more information on these
conjectures.

%%%%%%%%%%%%%%%%%%%%%%%%%%%%%%%%%%%%%%%%%%%%%%%%%%%%%%%%%%%%%%%%%%%%%%

%\input larmanr-bib.tex

\end{document}